\documentclass[12pt]{amsart}

\usepackage{amssymb}
\usepackage{latexsym}
\usepackage{verbatim}
\usepackage{fullpage}

\input {cyracc.def}
\font\ququ=cmr10 scaled \magstep1
\font\tencyr=wncyr8 scaled \magstep1

\def\rus{\tencyr\cyracc}

\newcommand{\re}[1]{\textrm  (\ref{#1})}

\renewenvironment{proof}
{\noindent {\sl Proof.}\quad }{\hfill $\square$ 
\vskip1.1ex\noindent }

\newenvironment{proof*}
{\noindent {\sl Proof.}\quad }{\hfill $\square$}

\renewcommand{\theequation}{\thesection.\arabic{equation}}
\renewcommand{\thesubsubsection}{\theequation.\arabic{subsubsection}}

\catcode`\@=\active
\catcode`\@=11
\def\@eqnnum{\hbox to
.01pt{}\rlap{\hskip-\displaywidth(\mathbf{\theequation})}}

\catcode`\@=12

\newenvironment{s}[1]
{ \vskip1.2ex \refstepcounter{equation}
\noindent {\bf \theequation\quad #1.} \begin{sl}}{\end{sl}
\vskip1.1ex\noindent }

\newenvironment{rem}[1]
{ \vskip1.2ex \refstepcounter{equation}
\noindent {\bf \theequation\enspace {#1}.} }{ \vskip1.1ex\noindent }

\newcommand {\ah}{{\frak a}}
\newcommand {\be}{{\frak b}}
\newcommand {\ce}{{\frak c}}

\newcommand {\g}{{\frak g}}

\newcommand {\te}{{\frak t}}

\newcommand {\z}{{\frak z}}

\newcommand {\esi}{\varepsilon}
\newcommand {\ap}{\alpha}

\newcommand {\tth}{\tilde{\theta}}

\newcommand {\ca}{{\mathcal A}}

\newcommand {\ck}{{\mathcal K}}
\newcommand {\N}{{\mathcal N}}
\newcommand {\co}{{\mathcal O}}

\newcommand {\md}{/\!\!/}

\newcommand {\isom}{\stackrel{\sim}{\longrightarrow}}
\newcommand {\ad}{{\mathrm{ad\,}}}
\newcommand {\Ad}{{\mathrm{Ad\,}}}

\newcommand {\Aut}{{\mathrm{Aut\,}}}

\newcommand {\ind}{{\mathrm{ind\,}}}
\newcommand {\Int}{{\mathrm{Int\,}}}

\newcommand {\Ker}{{\mathrm{Ker\,}}}

\newcommand {\rk}{{\mathrm{rk\,}}}
\newcommand {\spe}{{\mathrm{Spec\,}}}

\newcommand {\tri}{{\frak sl}_2}

\newcommand {\ov}{\overline}

\newcommand {\vno}[1]{\vskip#1 ex\noindent}

\newcommand {\beq}{\begin{equation}}
\newcommand {\eeq}{\end{equation}}

\newfam\Bbbfam\newfam\eufam%
\font\Bbbfont=msbm10 scaled 1200%
\font\olala=msam10 scaled 1200%
\font\frak=eufm10 scaled 1400%
\font\Bbbsmallfont=msbm8%
\textfont\Bbbfam=\Bbbfont\scriptfont\Bbbfam=\Bbbsmallfont
\font\euzw=eufm10 scaled 1200%
\font\euac=eufm7 scaled 1200%
\font\euacc=eufm7 scaled 1000%
\textfont\eufam=\euzw\scriptfont\eufam=\euac%
\scriptscriptfont\eufam=\euacc%
\def\frak{\fam\eufam}%
\def\Bbb{\fam\Bbbfam}%

\def\varnothing{\hbox {\Bbbfont\char'077}}
\def\square{\hbox {\olala\char"03}}
\def\bbk{\hbox {\Bbbfont\char'174}}

\begin{document}
\setlength{\parskip}{2pt plus 4pt minus 0pt}
\hfill {\scriptsize December 1, 2003} \vskip1ex
\vskip1ex

%
\title{On Invariant Theory of $\theta$-groups}
\author{{\sc Dmitri I. Panyushev}}
\maketitle
\vskip-2ex
\begin{center}
{\footnotesize 
\address{\it Independent University of Moscow,
Bol'shoi Vlasevskii per. 11 \\ 
121002 Moscow, \quad Russia \\ e-mail}: {\tt panyush@mccme.ru }\\
}
\end{center}

\hfill  {\it\small to Oksana}\hspace{10pt}

\section*{Introduction}
\vno{1}%
This paper is a contribution to Vinberg's theory of $\theta$-groups, or
in other words, to Invariant Theory of periodically graded semisimple
Lie algebras \cite{vi},\cite{vi79}. 
One of our main tools is Springer's theory of regular elements
of finite reflection groups \cite{tony}, with some recent complements by
Lehrer and Springer \cite{ls1},\,\cite{ls2}.
\\[.6ex]
The ground field $\bbk$ is algebraically closed and of 
characteristic zero.
Throughout, $G$ is a connected and simply connected semisimple algebraic 
group, $\g$ is its Lie algebra, and $\Phi$ is the Cartan--Killing form
on $\g$; $l=\rk\g$. 
\\
$\Int\g$ (resp. $\Aut\g$) is the group of inner (resp. all)
automorphisms of $\g$; $\N$ is the nilpotent cone in $\g$.
For $x\in\g$, $\z(x)$ is the centraliser of $x$ in $\g$.

Let $\g=\oplus_{i\in {\Bbb Z}_m}\g_i$ be a periodic grading of $\g$ and
$\theta$ the corresponding $m^{th}$ order automorphism of $\g$. 
Let $G_0$ denote the connected subgroup of $G$  with Lie 
algebra $\g_0$. 
Invariant Theory of $\theta$-groups deals with orbits and
invariants of $G_0$ acting on $\g_1$. Its main result is that there is a
subspace $\ce\subset\g_1$ and a finite reflection group $W(\ce,\theta)$
in $\ce$ (the {\it little Weyl group})
such that $\bbk[\g_1]^{G_0}\simeq \bbk[\ce]^{W(\ce,\theta)}$. 
We say that the grading is \textrm{N}-regular
(resp. S-regular) if $\g_1$ contains a regular nilpotent (resp. semisimple)
element of $\g$. The grading is locally free if there is $x\in \g_1$ such that
$\z(x)\cap\g_0=\{0\}$. The same terminology also applies to $\theta$.

In this paper, we obtain some structural results for
gradings with these properties and study interrelations of these properties.
Section~1 contains some preliminary material on $\theta$-groups and regular 
elements.
In Section~2, we begin with a dimension formula for semisimple $G_0$-orbits 
in $\g_1$. We also prove two ``uniqueness'' theorems. Recall that $\Int\g$ 
is the identity component of $\Aut\g$,  and it
operates on $\Aut\g$ via conjugations.
Given $m\in\Bbb N$, we prove that each connected component of $\Aut\g$
contains at most
one $\Int\g$-orbit consisting of automorphisms of order $m$ that are either
\textrm{N}-regular or S-regular and locally free.
In Section 3, we show that $\theta$-groups corresponding to
\textrm{N}-regular gradings enjoy a number of good properties.
Let $m_1,\dots,m_l$ be the exponents of $\g$ and let
$\{e,h,f\}$ be a regular $\tri$-triple such that $e\in\g_1$ 
and $f\in\g_{-1}$. For the purposes of this introduction, assume that
$\theta$ is inner. Set $k_i=\#\{j\mid m_j \equiv i \pmod{m} \}$
 ($i\in {\Bbb Z}_m$),
and let
$\zeta$ be a primitive $m^{th}$ root of unity.
We prove that 
(i) the eigenvalues of $\theta$ on $\z(e)$ are
$\zeta^{m_i}$ ($1\le i\le l$), 
(ii) $\dim\g_{i+1}-\dim\g_i=k_{-i-1}-k_i$ for all
$i\in {\Bbb Z}_m$, 
(iii) The restriction homomorphism $\bbk[\g]^G\to \bbk[\g_1]^{G_0}$ is
onto, (iv) the $G_0$-action on $\g_1$ admits a Kostant-Weierstrass (=\,KW)
section.
In the general case, the definition of the $k_i$'s becomes more involved,
see Eq.~\re{k_i},
but the above assertions (ii)--(iv) remains intact. 

In Section 4, it is shown that any locally free S-regular grading of $\g$ is 
\textrm{N}-regular. This implies that all such gradings admit a KW-section.
We also give a formula for dimension of all subspaces $\g_i$
in the S-regular case.
Another result is that $\dim\ce\le k_{-1}$ for any $\theta$-group.
We then show that the $G$-stable cone $\pi^{-1}\pi(\ce)\subset\g$ is 
a normal complete intersection. (Here $\pi:\g\to \g\md G$ is the 
quotient mapping.) In particular, if $\theta$
is S-regular or \textrm{N}-regular, 
then  $\ov{G{\cdot}\g_1}$ is a normal complete intersection.
A description of the defining ideal of $\ov{G{\cdot}\g_1}$ is also given.
This normality stuff relies on results of Richardson \cite{roger}.
It is curious to note that in case $m=2$ (i.e., $\theta$ is involutory) 
S-regularity is equivalent to \textrm{N}-regularity. But for $m>2$ neither of these
properties implies the other.   

Section 5 contains a description of the coexponents for little Weyl groups,
if $\theta$ is  both S- and \textrm{N}-regular. This is based on recent results of
Lehrer and Springer \cite{ls2}.

{\small {\bf Acknowledgements.} It is a pleasure for me to acknowledge the
benefit of conversations on $\theta$-groups that I had with E.B.\,Vinberg and
L.V.\,Antonyan.  This research was supported in part 
RFBI Grant 02--01--01041.
}


\section{Vinberg's $\theta$-groups and Springer's regular elements}


\noindent
Let $\theta$ be an automorphism of $\g$, of finite order $m$. The 
automorphism of $G$ induced by $\theta$ is also denoted by $\theta$.
Let $\zeta$ be a fixed primitive $m^{th}$ root of unity.
Then $\theta$ determines a periodic grading  
$\g=\displaystyle\bigoplus_{i\in {\Bbb Z}_m}\g_i$, where
$\g_i=\{x\in\g\mid \theta(x)=\zeta^i x\}$.
Whenever we want to indicate the dependence of the grading on
$\theta$, we shall endow `$\g_i$' with a suitable superscript.
If $M$ is a $\theta$-stable
subspace, then $M_i:=M\cap\g_i$.
Recall some standard facts on periodic gradings (see \cite[\S 1]{vi}):
\begin{itemize}
\item $\Phi(\g_i,\g_j)=0$ unless $i+j=0$;
\item $\Phi$ is non-degenerate on $\g_i\oplus\g_{-i}$ ($i\ne 0$) and
on $\g_0$. In particular,
$\g_0$ is a reductive algebraic Lie algebra and 
$\dim\g_i=\dim\g_{-i}$;
\item If $x\in\g_i$ and $x=x_s+ x_n$ is its Jordan decomposition, then
$x_s, x_n\in\g_i$.
\end{itemize}
Let $G_0$ be the connected subgroup of $G$ with Lie algebra $\g_0$. The 
restriction of the adjoint representation of $G$ to $G_0$ induces a 
representation $\rho_1$ of $G_0$ on $\g_1$. The linear group
$\rho_1(G_0)\subset GL(\g_1)$ is called a $\theta$-group. The theory of orbits
and invariants for $\theta$-groups, which generalizes that for the adjoint 
representation \cite{ko63} and for the isotropy representation of a symmetric
variety \cite{KR}, is developed by E.B.\,Vinberg in \cite{vi}.
\\[.8ex]
A {\it Cartan subspace\/} of $\g_1$ is a maximal commutative subspace
consisting of semisimple elements. Let $\ce\subset\g_1$ be a Cartan subspace.
Set $N(\ce)_0=\{g\in G_0\mid \Ad(g)\ce=\ce\}$ and
$Z(\ce)_0=\{g\in G_0\mid \Ad(g)x=x \textrm{ for all } x\in\ce\}$. The group
$N(\ce)_0/Z(\ce)_0$ is said to be the {\it little Weyl group\/}
of the graded Lie algebra, denoted $W(\ce,\theta)$. 

The following is a summary of main results in \cite{vi}.

\begin{s}{Theorem} \label{vinberg}
\begin{itemize}
\item[\sf (i)] \quad All Cartan subspaces in $\g_1$ are $G_0$-conjugate; 
\item[\sf (ii)] \quad $W=W(\ce,\theta)$ is a 
finite reflection group in $GL(\ce)$; 
\item[\sf (iii)] \quad  Let $x\in\g_1$. 
The orbit $G_0{\cdot}x$ is closed if and only if 
$G_0{\cdot}x\cap\ce\ne\varnothing$;
the closure $\ov{G_0{\cdot}x}$ contains the origin if and only if $x\in\N$; 
\item[\sf (iv)] \quad The restriction of polynomial functions
 $\bbk[\g_1]\to \bbk[\ce]$
induces an isomorphism $\bbk[\g_1]^{G_0}\isom \bbk[\ce]^W$
(a ``Chevalley-type'' theorem); 
\item[\sf (v)] \quad Each fibre of the quotient mapping $\pi_1:\g_1\to
\g_1\md G_0=\spe\bbk[\g_1]^{G_0}$ consists of finitely many $G_0$-orbits. 
The dimension of each fibre is equal to $\dim\g_1-\dim\ce$.
\end{itemize}
\end{s}%
Despite its maturity, the theory of $\theta$-groups still 
has a vexatious gap. A long-standing conjecture formulated in 
\cite[n.7]{cov}, to the effect that any $\theta$-group has an 
analogue of the section constructed by Kostant for the adjoint representation
(see \cite[n.\,4]{ko63}), is still open.
(In \cite{KR}, such a section was also constructed for
the isotropy representation of a symmetric variety. So that the problem
concerns the case $m\ge 3$.) 
Kostant's section for the adjoint representation
is an instance of a more general phenomenon in Invariant Theory, 
a so-called {\it Weierstrass section}.
The reader is referred to 
\cite[8.8]{vipo} for
the general definition of a Weierstrass section and a number of related
results. In my opinion, it is more natural to use term Kostant-Weierstrass 
sections, or {\it KW-sections} in the context of $\theta$-groups.
It was shown in \cite[Cor.\,5]{R=S} that a KW-section
exists whenever $\g_0$ is semisimple.  
Below, we will discuss some aspects of KW-sections in a more general 
situation.

If $W$ be a finite reflection group in a $\bbk$-vector space $V$, then
$v\in V$ is called regular if the stabiliser of $v$ in $W$ is
trivial. 
Let $\sigma$ be an element of finite order in $N_{GL(V)}(W)$.
Then $\sigma$ is called {\it regular \/}(in the sense of Springer)
if it has a regular eigenvector. The theory of such elements is developed
by Springer in \cite{tony}; for recent results, see \cite{ls1}.
Let $f_1,f_2,\dots,f_l$ be a set of algebraically independent homogeneous 
generators of $\bbk[V]^W$, $l=\dim V$. Set $d_i=\deg f_i$.
The $f_i$'s can be chosen so that $\sigma(f_i)=\esi_if_i$ ($1\le i\le l$),
with suitable roots of unity $\esi_i$.
Given a root of unity $\zeta$, we let $V(\sigma,\zeta)$ denote
the eigenspace of $\sigma$ corresponding to the eigenvalue $\zeta$.
The following is a sample of Springer's results, see Theorem~6.4
in \cite{tony}.

\begin{s}{Theorem}  \label{tony}
Suppose $V(\sigma,\zeta)$ contains a regular vector. Then
\begin{itemize}
\item[\sf (i)] $\dim V(\sigma,\zeta)=
\#\{j\mid 1\le j\le l,\ \esi_j\zeta^{d_j}=1\}$;
\item[\sf (ii)] The centraliser of $\sigma$ in $W$, denoted $W^\sigma$, 
is a reflection group in $V(\sigma,\zeta)$, and the restrictions 
$f_j\vert_{V(\sigma,\zeta)}$ with $\esi_j\zeta^{d_j}=1$ form a set of basic
invariants for $W^\sigma$;
\item[\sf (iii)] The eigenvalues of $\sigma$ in $V$ are
$\esi_i^{-1}\zeta^{-d_i+1}$, $i=1,\dots, l$;
\item[\sf (iv)] $\dim V(w\sigma,\zeta)\le \dim V(\sigma, \zeta)$ for all
$w\in W$.
If $\dim V(\sigma,\zeta)=\dim V(w\sigma,\zeta)$ for some
$w\in W$, then $\sigma$ and $w\sigma$ are conjugate by an element of $W$.
\end{itemize}
\end{s}%
We will apply Springer's theory in the context of $\theta$-groups, when
$V=\te$ is a $\theta$-stable
Cartan subalgebra of $\g$, $W=W(\te)$ is the Weyl group of $\te$,  
and $\sigma=\theta\vert_\te$. Obviously, such $\sigma$ normalises the
Weyl group.


\section{Miscellaneous results on periodic gradings}
\setcounter{equation}{0}

\begin{s}{Proposition}  \label{1/m} \begin{itemize}
\item[\sf (i)] If $x\in\g_1$ is semisimple, then
$\dim\g_k-\dim\z(x)_k$ does not depend on $k$. In particular,
$\dim [\g_0,x]=\frac{1}{m}\dim[\g,x]$;
\item[\sf (ii)] If $m=2$, then the relation 
$\dim [\g_0,x]=\frac{1}{2}\dim[\g,x]$ holds for all $x\in\g_1$.
\end{itemize}
\end{s}\begin{proof}
Consider the Kirillov form $\ck_x$ on $\g$. By definition,
$\ck_x(y,z)=\Phi(x,[y,z])$. 
From the invariance of $\Phi$ one readily deduces
that $\Ker\ck_x=\z(x)$. Since $x\in\g_1$, we have
$\ck_x(\g_i,\g_j)=0$ unless $i+j=-1$.
It follows that
\[
\dim\g_k-\dim\z(x)_k=\dim\g_{-k-1}-\dim\z(x)_{-k-1} \textrm{ for all }
k\in {\Bbb Z}_m \ .
\]
This already implies (ii). If $m$ is arbitrary, then one obtains
a good conclusion only for semisimple elements. Indeed,
if $\z(x)$ is reductive, then $\dim\z(x)_k=\dim\z(x)_{-k}$ for all
$k$. Hence $\dim\g_k-\dim\z(x)_k$ does not depend on $k$.
\end{proof}%
Part (ii) is due to Kostant and Rallis, see \cite[Prop.\,5]{KR}.
\\[.7ex]
{\bf Definition.} 
A periodic grading (or the corresponding automorphism)
of $\g$ is called {\rmfamily S}-{\it regular\/}, if
$\g_1$ contains a regular semisimple element of $\g$;
{\it \textup{\rmfamily N}-regular\/}, if
$\g_1$ contains a regular nilpotent element of $\g$.
It is called {\it locally free\/}, if there exists $x\in\g_1$ such that
$\z(x)_0=\{0\}$.

Our aim is to prove two conjugacy theorems for periodic gradings.
  
\begin{s}{Theorem}  \label{conj-ss}
Let $\theta_1$, $\theta_2$ be automorphisms of $\g$ having the same order
and lying in the same connected component of $\Aut\g$.
Suppose the corresponding periodic gradings are \textup{\rmfamily S}-regular 
and locally-free.
Then $\theta_1$, $\theta_2$ are conjugated by means of an element of\/
$\Int\g$.
\end{s}\begin{proof}
Let $\g=
\bigoplus_{i=0}^{m-1}\g_i^{(1)}$ and 
and $\g=
\bigoplus_{i=0}^{m-1}\g_i^{(2)}$ be the corresponding 
gradings. It follows from the hypotheses that 
$\g_1^{(i)}$ ($i=1,2$) contains a 
regular semisimple element $x_i$ such that $\z(x_i)_0=\{0\}$.
Since all Cartan subalgebras are conjugate with respect to $\Int\g$, 
we may assume 
that $\z(x_1)=\z(x_2)=:\te$. Set 
$\sigma_i=\theta_i\vert_\te\in\Aut(\te)$.
Since $\theta_1\theta_2^{-1}$ is inner, we see that 
$\sigma_1=w\sigma_2 $ for some $w\in W(\te)$, where $W(\te)$ is the Weyl group
of the pair $(\g,\te)$. Thus, $\sigma_2$ and $w\sigma_2$ are two elements on
$GL(\te)$ having a regular eigenvector, with the same eigenvalue.
By \cite[6.4(iv)]{tony}, it then follows that $\sigma_1$ and $\sigma_2$ are 
$W(\te)$-conjugate. Thus, we may assume that $\theta_1\vert_\te=
\theta_2\vert_\te$. Then $\theta_2\theta_1^{-1}=\Ad (s)$ for some $s\in T
=\exp(\te)$.
By assumption, $\te_0=\{0\}$. Therefore $\theta_2\vert_T$ has finitely many
fixed points. Hence the mapping 
$(g\in T) \mapsto (\theta_2(g)g^{-1}\in T)$ is onto, 
and there exists $t\in T$
such that $s=\theta_2(t)t^{-1}$. Then
\[
\theta_1=\Ad (s)^{-1}{\cdot}\theta_2=\Ad (t){\cdot}\Ad(\theta_2(t^{-1}))
{\cdot}\theta_2=
\Ad (t){\cdot}\theta_2{\cdot}\Ad (t^{-1}) \ ,
\]
and we are done.
\end{proof}%
Recall that $x\in\N$ is called {\it semiregular\/}, if any semisimple 
element of the centraliser $Z_G(x)$ belong to the centre of $G$. The corresponding
orbit and $\tri$-triple are also called semiregular.
The semiregular $\tri$-triples in simple Lie algebras were classified 
by E.B.~Dynkin in 1952.

\begin{s}{Theorem}  \label{conj-nilp}
Let $\theta'$, $\theta''$ be automorphisms of $\g$ 
having the same order and 
lying in the same connected component of $\Aut\g$.
Suppose there exists a semiregular nilpotent orbit $\co\in\g$ such that
$\co\cap\g'_1\ne\varnothing$ and $\co\cap\g''_1\ne\varnothing$.
Then $\theta'$, $\theta''$ are conjugated by means of an element of\/
$\Int\g$.
\end{s}\begin{proof*}
In case $\co$ is the regular nilpotent orbit, a proof is given in
\cite{leva2}. It goes through in our slightly more general setting. 
For convenience of the reader, we give it here.

Let $e'\in\g'_1\cap\co$ and $e''\in\g_1''\cap\co$.
According to \cite[sect.\,2]{vi79}, there exist $\tri$-triples
$\{e',h',f'\},\ \{e'',h'',f''\}$ such that
$h'\in\g_0',\, h''\in\g_0''$, $f'\in\g_{-1}',\, f''\in\g_{-1}''$. 
By the conjugacy theorem for $\tri$-triples,
there exists $\tau\in\Int\g$ such that $\tau(e')=e'',\,\tau(h')=h'',\,
\tau(f')=f''$. Then ${\theta'}^{-1}\tau^{-1}\theta''\tau$ is inner, and it
takes the  triple $\{e',h',f'\}$ to itself. Since the centraliser of
a semiregular $\tri$-triple in $\Int\g$ is trivial, we obtain
$\theta'=\tau^{-1}\theta''\tau$.
\end{proof*}%

\section{\textrm{N}-regular periodic gradings } \label{section-R}
\setcounter{equation}{0}

\noindent
Let $\co^{reg}$ be the
regular nilpotent orbit in $\g$. Recall that a periodic grading 
(or automorphism) of $\g$ is \textrm{N}-regular, if 
$\co^{reg}\cap\g_1\ne\varnothing$.
Since $\co^{reg}$ is semiregular, Theorem~\ref{conj-nilp} says
that any connected component of $\Aut\g$ 
contains at most one $\Int\g$-orbit of \textrm{N}-regular automorphisms of a 
prescribed order.
To give a detailed description of the \textrm{N}-regular periodic gradings,
some preparatory work is needed.

For any $\gamma\in\Gamma(\g):=\Aut\g/\Int\g$, 
let $C_\gamma$ denote the corresponding
connected component of $\Aut\g$. The {\it index\/} of 
(any element of) $C_\gamma$ is the order of $\gamma$ in $\Gamma(\g)$.
The index of $\mu\in\Aut\g$ is denoted by $\ind\mu$. Thus,  
$\textrm{ord}\gamma=\ind C_\gamma=\ind\mu$ for any $\mu\in C_\gamma$.

Since $\Int\g\simeq G/\{centre\}$, the group $\Gamma(\g)$ acts on
$\bbk[\g]^G$ (or on $\g\md G=\spe\bbk[\g]^G$).
Let $\mu\in 
\Aut\g$ be arbitrary. 
 Denote by $\ov{\mu}$ the corresponding 
(finite order) automorphism of $\g\md G$. 

\begin{s}{Lemma} \label{order}
The action of $\Gamma(\g)$ on $\g\md G$ is effective. In other words,
the order of\/ $\ov{\mu}$ equals $\ind\mu$.
\end{s}\begin{proof}
It is clear that the order of $\ov{\mu}$ divides  $\ind\mu$. 
To prove the converse, we have to show that if $\ov{\mu}$ is trivial, then
$\mu$ is inner. Without loss of generality, one may assume that $\mu$ is a
semisimple automorphism. Then,
by a result of Steinberg \cite[Thm.\,7.5]{endo}, 
there is a Borel subalgebra
$\be\subset\g$ and a Cartan subalgebra $\te\subset\be$ such that 
$\mu(\be)=\be$ and $\mu(\te)=\te$. Let $W(\te)$ be the Weyl group of $\te$.
Since $\mu$ acts trivially on $\bbk[\g]^G\simeq \bbk[\te]^{W(\te)}$,
the restriction of $\mu$ to $\te$ is given by an element of $W(\te)$.
On the other hand, the relation $\mu(\be)=\be$ shows that $\mu\vert_\te$
permutes somehow the simple roots corresponding to $\be$.
It follows that $\mu$ acts trivially on $\te$ and
therefore $\mu$ is inner.
\end{proof}%
In the following theorems, we describe \textrm{N}-regular periodic 
gradings  and give some relations for eigenvalues and eigenspaces of
$\theta$.
\begin{s}{Theorem {\ququ (Antonyan)}}  \label{Antonyan}
Fix $m\in\Bbb N$, and consider a connected component $C_\gamma\subset\Aut\g$.
Then \\[.7ex]
\parbox{\textwidth}{\hfil
$\left\{
\parbox{165pt}{$C_\gamma$ contains an \textup{\rmfamily N}-regular automorphism 
of order $m$
}\right\}$ \quad $\Longleftrightarrow$ \quad
\parbox{150pt}{$\ind C_\gamma$ divides $m$.}
\hfil}
\vskip.5ex\noindent
In other words, if a connected component of $\Aut\g$ contains elements
of order $m$, then it contains an \textup{\rmfamily N}-regular automorphism 
of order $m$. 
\end{s}\begin{proof}
If $\ind C_\gamma$ does not divide $m$, then $C_\gamma$ does not contain
automorphisms of order $m$. 
To prove the converse, we first fix a Borel subalgebra $\be\subset\g$
and a Cartan subalgebra $\te\subset\be$. Let $\Delta$ be the root system
of $(\g,\te)$. Let $\Pi=\{\ap_1,\dots,\ap_l\}$ be the set of simple roots 
such that the roots of $\be$ are positive. For each $\ap_i\in\Pi$, let 
$e_i$ be a nonzero root vector.
Recall that the finite group $\Gamma(\g)$ is isomorphic to
the symmetry group of the Dynkin diagram of $\g$ \cite[4.4]{vion}.
This means that each $C_\gamma$ contains an automorphism 
$\theta_\gamma$ such that 
$\theta_\gamma(\te)=\te$ and $\theta_\gamma(e_i)=c_ie_{\bar\gamma(i)}$
, $i=1,\dots,l$,
where $\bar\gamma$ is a permutation on $\{1,\dots,l\}$ and 
$c_i\in\bbk\setminus\{0\}$. The permutation $\bar\gamma$ represents 
an automorphism of the Dynkin diagram of $\g$ and the order of
$\bar\gamma$ equals $\ind\mu$. Conjugating $\theta_\gamma$ by 
$\Ad(t)$ for a suitable $t\in T$, we can obtain arbitrary coefficients
$c_i$. Therefore we may assume without loss of generality that
$c_1=\dots=c_l=\zeta$. Then $\theta=\theta_\gamma$ is \textrm{N}-regular, and
of order $m$. Indeed, 
$e_1+\dots +e_l$ is a regular nilpotent
element lying in $\g_1^{(\theta)}$. Next, $\theta^m$ is inner and
$\theta^m(e_i)=e_i$ for all $i$. Hence $\theta^m=id_\g$.
\end{proof}%
Let $F_1,F_2,\dots,F_l$ be homogeneous algebraically independent 
generators of $\bbk[\g]^G$, 
$\deg F_i=d_i$. Set $m_i=d_i-1$. The numbers $m_1\dots,m_l$ are called the
{\it exponents\/} of $\g$. 
Given $\theta\in C_\gamma\subset\Aut\g$, we may choose the $F_i$'s so that
$\theta(F_i)=\esi_i F_i$ ($i=1,\dots,l$) for some roots of unity $\esi_i$.
We shall say that the $\esi_i$'s are the {\it factors\/} of $\theta$. 
Note that the multiset $\{\esi_1,\dots,\esi_l\}$ depends only on the 
connected component of $\Aut\g$, containing $\theta$.
If $\te$ is an arbitrary
Cartan subalgebra, then $f_i=F_i\vert_\te$ ($1\le i\le l$) are algebraically
independent generators for $\bbk[\te]^{W(\te)}$ and $\sigma(f_i)=\esi_i f_i$,
where $\sigma=\theta\vert_\te$. So, the $\esi_i$'s are also factors in the
sense of Springer \cite[\S\,6]{tony}.
Given $m\in\Bbb N$, we shall exploit two
sequences indexed by elements of ${\Bbb Z}_m$.
Set 
\begin{equation}  \label{k_i}
k_{i}:=\#\{j\mid 1\le  j\le l,\ \zeta^{m_j-i}\esi_j=1 \}
\quad \textrm{and} \quad
l_{i}:=\#\{j\mid 1\le  j\le l,\  \zeta^{m_j-i}\esi_j^{-1}=1 \} \ .
\end{equation}
In this way, we obtain the numbers 
satisfying the relation $\sum_{i\in {\Bbb Z}_m} k_i=
\sum_{i\in {\Bbb Z}_m} l_i=l$.
Given $\theta\in\Aut\g$ \ ($\theta^m=id_\g$), an 
$\tri$-triple $\{e,h,f\}$ is said to be $\theta$-{\it adapted}, if 
$\theta(e)=\zeta e$, $\theta(h)=h$, and
$\theta(f)=\zeta^{-1}f$. By an extension of the Morozov-Jacobson theorem
\cite[\S\,2]{vi79},
any nilpotent $e\in\g_1$ can be included in a $\theta$-adapted $\tri$-triple.
Recall that an $\tri$-triple is called {\it regular\/}, if $e\in \co^{reg}$.

\begin{s}{Theorem}  \label{numerology}
Suppose $\theta\in\Aut\g$ is \textup{\rmfamily N}-regular and of order $m$.
Let $\{\esi_1,\dots,\esi_l\}$ be the factors of $\theta$ and
$\{e,h,f\}$  a $\theta$-adapted regular $\tri$-triple.
Then 
\begin{itemize}
\item[\sf (i)]\ The eigenvalues of $\theta$ on $\z(e)$ are equal to
$\zeta^{m_i}\esi_i^{-1}$, $i=1,\dots,l$;
\item[\sf (ii)]\ The eigenvalues of $\theta$ on $\z(h)$ (resp.
$\z(f)$) are equal to $\esi_i^{-1}$ (resp. $\zeta^{-m_i}\esi_i^{-1}$), 
$i=1,\dots,l$;
\item[\sf (iii)]\ $k_i=l_i$ for all $i\in {\Bbb Z}_m$.
\item[\sf (iv)]\ $ \dim\g_{i+1}-\dim\g_i=k_{-i-1}-k_i$.
\item[\sf (v)]\ The dimension of a Cartan subspace of $\g_1$ equals $k_{-1}$.
\end{itemize}
\end{s}\begin{proof}
Set $\ah=\bbk e+\bbk h+ \bbk f\simeq\tri$. 
Let $R(n)$ denote the irreducible
$\ah$-module of dimension $n+1$. It is well known that
as $\ah$-module $\g$ is isomorphic to $\oplus_{i=1}^l R(2m_i)$. 

(i) Let $\partial_f$ be the derivation of $\bbk[\g]$
determined by $f$. Then $q_i=(\partial_f)^{m_i}F_i$ is a linear
form on $\g$. By \cite[Theorem~6]{ko63}, there exists
a basis $x_1,\dots,x_l$ for $\z(e)$ such that $[h,x_i]=2m_ix_i$ and
$q_i(x_j)=\delta_{ij}$.
Since $\theta(q_i)=\zeta^{-m_i}\esi_i$, the latter implies that
$\theta(x_i)=\zeta^{m_i}\esi_i^{-1}$.

(ii) Since  $x_i$ is a highest weight vector in $R(2m_i)$, the
vectors $h_i:=(\ad f)^{m_i}x_i$ and
$y_i:=(\ad f)^{2m_i}x_i$ $(i=1,\dots,l)$
form a basis for $\z(h)$ and $\z(f)$, respectively.
Obviously, $\theta(h_i)=\esi_i^{-1}$ and 
$\theta(y_i)=\zeta^{-m_i}\esi_i^{-1}$. 

(Another proof for $\z(h)$ can be derived from \cite{tony}.
Since $h$ is regular semisimple, $\z(h)$ is a $\theta$-stable 
Cartan subalgebra. Set $\sigma=\theta\vert_{\z(h)}$.
Then $\sigma$ is a regular element in the sense
of Springer. Indeed, $\sigma$ normalizes the Weyl group
$N_G(\z(h))/Z_G(\z(h))$ and $h$ is a regular eigenvector of $\sigma$.
By \cite[6.5(i)]{tony}, the eigenvalues of $\sigma$ on $\z(h)$
are equal to $\esi_i^{-1}$, $i=1,\dots,l$.)

(iii) It follows from (i) and (ii) that $\dim\z(e)_i=l_i$
and $\dim\z(f)_{-i}=k_i$. 
It remains to observe that
$\dim\z(e)_i=\dim\z(f)_{-i}$, since $\Phi$ is $\theta$-invariant and
yields a nondegenerate pairing
between $\z(e)$ and $\z(f)$.


(iv) Since $x_i$ is a highest
weight vector in $R(m_i)$, it follows from (i) that
the eigenvalues of $\theta$ in $\g$ are
\[
  \{ \zeta^i\esi_j^{-1}\mid j=1,2,\dots,l;\  i=m_j,m_j-1,\dots,-m_j\} \ .
\]
So, the problem of computing the required differences becomes purely 
combinatorial.
Notice that each $\esi_i$ is an $m^{th}$ root of unity, so that each eigenvalue
is a power of $\zeta$. Let us calculate separately the contribution of each
submodule $R(2m_j)$ to the difference 
$D_i:=\dim\g_{i+1}-\dim\g_i$.
Usually, two consecutive eigenvectors with eigenvalues 
$\zeta^{i+1}$ and $\zeta^i$ occur together in $R(2m_j)$;
i.e., this has no affect on the difference in question.
The exceptions can only occur near the eigenvalues of the highest  and
the lowest weight vectors in $R(2m_j)$. Namely, if 
$\zeta^{m_j}\esi_j^{-1}=\zeta^i$, one gains contribution $-1$ to $D_i$; 
if $\zeta^{-m_j}\esi_j^{-1}=\zeta^{i+1}$, one gains contribution 
$+1$ to $D_i$. Thus, taking the sum over all irreducible $\ah$-submodules
yields
\[
D_i=\#\{ j\mid \zeta^{-m_j-i-1}\esi_j^{-1}=1\} -
  \#\{ j\mid \zeta^{m_j-i}\esi_j^{-1}=1\}=k_{-i-1}-l_i=
  k_{-i-1}-k_i \ .
\]
(v) 
By parts (i) and (iv), we have 
$\dim\z(e)_0=\#\{j\mid\zeta^{m_j}\esi_j^{-1}=1\}=l_0=k_0$ 
and $\dim\g_{1}-\dim\g_0=k_{-1}-k_0$. 
Therefore $\dim G_0{\cdot}e=
\dim\g_0- k_0$. Since $G{\cdot}e$ is open and dense in $\N$, $G_0{\cdot}e$ is a
nilpotent orbit in $\g_1$ of maximal dimension. By Theorem~\ref{vinberg}(v),
$\dim G_0{\cdot}e$ is also maximal among dimensions of
{\sl all\/} $G_0$-orbits in $\g_1$.
Thus,
\[
\dim\ce=\dim\g_1\md G_0=\dim\g_1-\dim G_0{\cdot}e=k_{-1} \ .
\]\end{proof}%

In the next claim, we regard $\{0,1,\dots,m-1\}$ as a
set of representatives for ${\Bbb Z}_m$.
\begin{s}{Corollary}
$\dim\g_0=\displaystyle
\frac{1}{m}\bigl(\dim\g+\sum_{i=0}^{m-1}(m-1-2i)k_i\bigr)$.
\end{s}\begin{proof}
Write the relations of Theorem~\ref{numerology}(iv) in the form
$\dim\g_{i+1}-\dim\g_i=k_{m-i-1}-k_i$,
$0\le i\le m-1$. Together
with the equality $\sum_{i=0}^{m-1}\dim\g_i=\dim\g$, these form a system of 
$m$ linear equations with $m$ indeterminates $\{\dim\g_i\}$.
\end{proof}%
Utility of \textrm{N}-regularity is explained by the fact that this 
allows us describe the algebra of invariants $\bbk[\g_1]^{G_0}$
and guarantee the
existence of a KW-section. Let us briefly recall the last subject.
An affine subspace $\ca\subset \g_1$ is called a KW-section if the restriction
of $\pi_1:  \g_1\to \g_1\md G_0$ to $\ca$ is an isomorphism.
By Theorem~\ref{vinberg}(iii), such an $\ca$ contains a unique nilpotent 
element. So that $\ca$ is of the form $v+L$, where $\{v\}=\ca\cap\N$ and $L$
is a linear subspace of $\g_1$.


\begin{s}{Theorem}  \label{onto+KW}
Suppose $\theta$ is \textup{\rmfamily N}-regular and of order $m$. 
Let $\{e,h,f\}$  a $\theta$-adapted regular $\tri$-triple. Then
\begin{itemize} 
\item[\sf (i)]\ 
The restriction homomorphism $\bbk[\g]^G\to \bbk[\g_1]^{G_0}$ 
is onto. Moreover, $\bbk[\g_1]^{G_0}$ is freely generated by the restriction
to $\g_1$ of all basic invariants $F_j$ such that $\zeta^{m_j+1}\esi_j=1$.
\item[\sf (ii)]\ $e+\z(f)_1$ is a KW-section in $\g_1$.
\end{itemize}
\end{s}\begin{proof} 
(i) Choose a numbering of basic invariants so that the relation
$\zeta^{m_i+1}\esi_i=1$ holds precisely for $i\le a$. 
Observe that $a=k_{-1}$.
It is immediate that $F_i$ vanishes on $\g_1$ unless $i\le k_{-1}$. For,
\[
\esi_iF_i(x)=(\theta (F_i))(x)=F_i(\theta^{-1}(x))=F_i(\zeta^{-1} x)
=\zeta^{-m_i-1}F_i(x)\ \textrm{ for all }\ x\in\g_1 \ .
\]
(Recall that $d_i=m_i+1$.) Our aim is to show that
$\bar F_i:=F_i\vert_{\g_1}$ ($1\le i\le k_{-1}$) generate $\bbk[\g_1]^{G_0}$.
(By Theorem~\ref{vinberg}, $\bbk[\g_1]^{G_0}$ is a polynomial algebra in
$\dim\ce$ variables, i.e., in our case in $k_{-1}$ variables.)
A standard fact of $\tri$-theory says that $\z(f)\oplus [\g,e]=\g$.  By a
famous result of Kostant \cite{ko63}, $(dF_i)_e$ are linearly independent
as elements of $\g^\ast$ and their images in $\z(f)^\ast$ form a basis 
for $\z(f)^\ast$. Therefore, restricting the differentials of basic
invariants to $\z(f)_1$, one obtains a basis for $\z(f)_1^\ast$.
The preceding exposition shows 
$(d\bar F_i)_e=(d F_i)_e\vert_{\g_1}=0$ unless $1\le i\le k_{-1}$. On the
other hand, it follows from Theorem~\ref{numerology}(i) that
$\dim\z(f)_1=k_{-1}$. 
Hence $(d\bar F_i)_e$ ($1\le i\le k_{-1}$) are linearly independent 
and $\bar F_i$ are algebraically independent.
Furthermore, the linear independence of differentials
implies that each $\bar F_i$ is a member of minimal generating
system for $\bbk[\g_1]^{G_0}$, since
$e$ lies in the zero locus of all homogeneous $G_0$-invariants of positive
degree. This completes the proof.

(ii) It is a standard consequence of
the fact that $\{\bar F_i\}$ generate $\bbk[\g_1]^{G_0}$ and
$(d\bar F_i)_e$ ($1\le i \le k_{-1}$) are linearly independent, 
see e.g. \cite[\S\,3]{R=S}. 
\end{proof}%
%
%
{\it Remark.} If $\theta$ is inner, then $\esi_i=1$ for all $i$, and
the previous exposition simplifies considerably. In this case, we also
have $k_i=\#\{ j\mid m_j \equiv i\pmod{m}\}$.

\begin{s}{Corollary}
Suppose $\theta$ is inner and \textup{\rmfamily N}-regular. Then $\bbk[\g_1]^{G_0}$
is freely generated by  those $F_i\vert_{\g_1}$ whose degree is divisible my 
$m$.
\end{s}%
\vskip-1ex

\section{Applications and examples}
\setcounter{equation}{0}

\noindent
In this section, we demonstrate some applications of Springer's theory
of regular elements to $\theta$-groups.
Maintain the notation of the previous section. In particular,
to any $\theta\in\Aut\g$, of order $m$,
we associate the factors $\esi_i$ $(i=1,\dots,l)$, which depend only
on
the connected component of $\Aut\g$ that contains $\theta$, and then
the numbers $k_i=k_i(\theta,m)$ $(i\in {\Bbb Z}_m)$, 
which are defined by Eq.~\re{k_i}.

\begin{s}{Lemma}   \label{k-1}
Let $\theta\in\Aut\g$ be any automorphism of order $m$.
Then $\dim\g_1\md G_0\le k_{-1}$.
\end{s}\begin{proof}
Take a Cartan subspace $\ce\subset\g_1$.  Let $\te$ be any $\theta$-stable
Cartan subalgebra of $\g$ containing $\ce$. 
Because of the maximality of $\ce$, we have $\te_1=\ce$.  Set $\sigma=
\theta\vert_{\te}$. Then $\sigma\in N_{GL(\te)}(W(\te))$ and
$\ce=\te(\sigma,\zeta)$. By \cite[6.2(i)]{tony}, one has
$\dim\ce\le k_{-1}$ (our $k_{-1}$ is $a(d,\sigma)$ in Springer's paper).
\end{proof}%
It is shown in the previous section that the $k_i$'s play a significant
r\^ole in the context of \textrm{N}-regular gradings. Now we show that these
numbers also relevant to S-regular gradings.
\begin{s}{Theorem}  \label{main}
Let $\g=\displaystyle\bigoplus_{i\in {\Bbb Z}_m}\g_i$ be an \textup{\rmfamily S}-regular
grading and $\theta$ the corresponding automorphism. 
Then 
\begin{itemize}
\item[\sf (i)] $k_i=k_{-i}$ and $\dim\g_1\md G_0=k_{-1}$.
\item[\sf (ii)] $\dim\g_0-k_0=(\dim\g-l)/m$, and $\dim\g_i+k_0=\dim\g_0+k_i$.
\item[\sf (iii)] If $\theta$ is also locally free, then
$\co^{reg}\cap\g_1\ne\varnothing$, i.e., $\theta$ is\/ \textup{\rmfamily N}-regular.
\end{itemize}
\end{s}\begin{proof}
1.
Let $x\in\g_1$ be a regular semisimple element.
Set $\te=\z(x)$. It is
a $\theta$-stable Cartan subalgebra of $\g$. Set $\sigma:=\theta\vert_\te$.
Because $\sigma$ originates from an automorphism of
$\g$, it normalizes $W(\te)$, the Weyl group of $\te$.
Furthermore $x$ is a regular eigenvector of $\sigma$ whose eigenvalue is
$\zeta$. Thus, $\sigma$ is a regular element of $GL(\te)$
in the sense of Springer, and we conclude from \cite[6.4(v)]{tony} that
the eigenvalues of $\sigma$ are equal to $\zeta^{-m_i}\esi_i^{-1}$
($i=1,\dots,l$).
It follows that $\dim\te_i=k_{-i}$. Since $\Phi\vert_\te$ is nondegenerate
and $\sigma$-stable, $k_{-i}=k_i$. It is also clear that $\te_1$ is a
Cartan subspace of $\g_1$.

2. These two relations follow from Proposition~\ref{1/m}(i) 
applied to $x$.

3. As the grading is locally-free, $k_0=\dim\z(x)_0=0$.
Let $\tth$ be an \textrm{N}-regular automorphism of order $m$  that
lies in the {\sl same\/} connected component of $\Aut\g$ as $\theta$
(cf. Theorem~\ref{Antonyan}).
Let $\g=\oplus_i\tilde\g_i$ be the corresponding grading.
Let $\ce\subset\tilde\g_1$ be a Cartan subspace. By Theorem~\ref{numerology}(v),
$\dim\ce=k_{-1}$. Let $\tilde\te$ be any $\tth$-stable Cartan subalgebra 
containing $\ce$. Then, by the definition of a Cartan subspace, we have
$\ce=\tilde\te_1$. Conjugating $\tth$ by a suitable inner automorphism,
we may assume that $\te=\tilde\te$. Set $\tilde\sigma=\tth\vert_{\te}$.
Since $\theta\tth^{-1}$ is inner by the construction, 
$\sigma\tilde\sigma^{-1}\in W(\te)$.
Thus, we have the following:

{\it $\sigma$ has finite order,
$\sigma W(\te) \sigma^{-1}=W(\te)$, $\sigma=w\tilde\sigma$ for some 
$w\in W(\te)$,
and $\dim\te_1=\dim\tilde\te_1$}. 
\\
Since $\te_1$ contains a regular vector,
Theorem~6.4(ii),(iv) from \cite{tony} applies. 
It asserts that $\sigma$ and $\tilde\sigma$
are conjugate by an element of $W(\te)$. 
It follows that $\dim\te_i=\dim\tilde\te_i$ for all $i$, and
$\tilde\te_1$ contains a regular vector, too.
Thus, $\tth$ is S-regular and locally free as well.
Finally, applying Theorem~\ref{conj-ss} to $\theta$ and $\tth$, we conclude 
that these two are conjugate by an element of $\Int\g$.
Hence $\theta$ is also \textrm{N}-regular. 
\end{proof}%
{\it Remark.} If $\theta$ is not assumed to be locally free, then
part (iii) can be false, see example below.

Combining Theorem~\ref{onto+KW}(ii) and Theorem~\ref{main}(iii), we obtain

\begin{s}{Corollary}
If $\theta$ is\/ \textup{\rmfamily S}-regular and locally free, then the corresponding 
$\theta$-group admits a KW-section.
\end{s}%
Examples show that $\N\cap\g_1$, the null-fibre of $\pi_1$,
is often reducible. Any KW-section, if it exists, must meet one of the
irreducible components of $\N\cap\g_1$.  It turns out,
however, that some components are `good' and some are `bad' in this sense.
It may happen that there is only one irreducible component that can be
used for constructing a KW-section. It is worth noting in this regard that,
in case $\theta$ is involutory, 
all irreducible components of $\N\cap\g_1$ are `good', see
\cite[Theorem\,6]{KR}.
 
\begin{rem}{Example}  \label{primer}
Let $\g$ be a simple Lie algebra of type ${\Bbb E}_6$.
Consider two inner automorphisms $\theta_1,\theta_2$ of $\g$ 
that are defined by the following Kac's diagrams: \\
\centerline{
\raisebox{20\unitlength}{$\theta_1$:}\quad \begin{picture}(70,55)(0,0)
\setlength{\unitlength}{0.017in}
\multiput(35,3)(0,30){2}{\circle*{5}} 
\put(35,18){\circle{5}}
\multiput(5,33)(15,0){5}{\circle{5}}
\multiput(8,33)(15,0){4}{\line(1,0){9}}\multiput(35,6)(0,15){2}{\line(0,1){9}}
\end{picture}
\hspace{2cm}
\raisebox{20\unitlength}{$\theta_2$:}\quad \begin{picture}(70,55)(0,0)
\setlength{\unitlength}{0.017in}
\multiput(35,3)(0,30){2}{\circle{5}} 
\put(35,18){\circle*{5}} 
\multiput(5,33)(60,0){2}{\circle*{5}}
\multiput(20,33)(30,0){2}{\circle{5}}
\multiput(8,33)(15,0){4}{\line(1,0){9}}\multiput(35,6)(0,15){2}{\line(0,1){9}}
\end{picture}
}

\noindent
The reader is referred to \cite[4.7]{vion} or \cite[\S\,8]{vi} for a thorough
treatment of Kac's diagrams of periodic automorphisms.
Here we give only partial explanations:
\begin{itemize}
\item (The conjugacy class of) a periodic {\it inner\/} automorphism of a 
simple Lie algebra
$\g$ is represented by the corresponding {\it affine\/} Dynkin diagram,
with white and black nodes.
\item The semisimple part of $\g_0$ is given by the subdiagram consisting
of white nodes.
\item Dimension of the centre of $\g_0$ equals the number of black nodes
minus 1.
\item The order of $\theta$ is equal to the sum of those coefficients of
the affine Dynkin diagram that correspond to the black nodes.
\item Each black node represents an irreducible $\g_0$-submodule of $\g_1$,
so that the number of black nodes is equal to the number of 
irreducible summands of $\g_1$.
(We do not give here a general recipe for describing the
$\g_0$-module $\g_1$.)
\end{itemize}
It follows that both automorphisms under consideration have order 4,
$G_0^{(1)}=A_2\times A_2\times A_1\times\bbk^*$, and
$G_0^{(2)}=A_3\times A_1\times(\bbk^*)^2$. 

The $G_0^{(1)}$-module
$\g_1^{(1)}$ has 2 summands: tensor product of simplest representations
of all simple factors (dimension 18) plus 2-dimensional representation 
of $A_1$.
The weights of $\bbk^*$ on these summands, say $\mu_1$ and $\mu_2$, 
satisfy the relation $3\mu_1+\mu_2=0$ (in the additive notation).

The $G_0^{(2)}$-module
$\g_1^{(2)}$ has 3 summands: simplest representation of $A_3$ plus 
its dual plus
tensor product of the second fundamental representation of $A_3$ and
the simplest representation of $A_1$.
The weights of $(\bbk^*)^2$ on these summands, say $\mu_i$ ($i=1,2,3$),
satisfy the relation $\mu_1+\mu_2+2\mu_3=0$ (in the additive notation).

We have $\dim\g_0^{(i)}=\dim\g_1^{(i)}=20$, $i=1,2$. 
A direct computation shows in both cases that
the action $G_0^{(i)}:\g_1^{(i)}$ is stable, and stabiliser in general position is
a 2-dimensional torus. Hence $\dim G_0^{(i)}{\cdot}x_i=18$ for a generic (semisimple)
$x_i\in\g_1^{(i)}$ and therefore $\dim G{\cdot}x_i=72$ by 
Proposition~\ref{1/m}. Notice that $72=\dim\g-\rk\g$.
Thus, both $\theta_1$ and $\theta_2$ are S-regular but not locally free. 
Clearly, these are not conjugate. This proves that the assumption of being
locally free cannot be dropped in Theorem~\ref{conj-ss}.
\\[.6ex]
It is not hard to compute directly that the degrees of basic 
$G_0^{(i)}$-invariants
are equal to 8,\,12 for $\theta_1$  and 4,\,8 for $\theta_2$.
As the degrees of ${\Bbb E}_6$ are $2,\,5,\,6,\,8,\,9$ and $12$, we see that
the restriction $\bbk[\g]^G\to \bbk[\g_1^{(2)}]^{G_0^{(2)}}$ is not onto.
Hence $\theta_2$ is not \textrm{N}-regular. This proves that the assumption 
of being locally free cannot be dropped in Theorem~\ref{main}(ii).
By the way, $\theta_1$ is \textrm{N}-regular.
\end{rem}%

\noindent
Let $\theta$ be an arbitrary periodic automorphism and let
$G_0:\g_1$ be the corresponding $\theta$-group, with a Cartan subspace
$\ce\subset\g_1$ and the little Weyl group $W(\ce,\theta)$. 
The isomorphism \ref{vinberg}(iv) means that
$G_0{\cdot}x\cap\ce=W(\ce,\theta){\cdot}x$ for all $x\in\ce$.
It is not however always  true that $G{\cdot}x\cap\ce=G_0{\cdot}x\cap\ce$ for 
all $x\in\ce$.
A similar phenomenon can be seen on the level of Weyl groups, as follows. 
Let $\te$ be a $\theta$-stable Cartan subalgebra such that $\te_1=\ce$.
Write $W$ for the Weyl group $N_G(\te)/Z_G(\te)$. 
Set $W_1=N_W(\ce)/Z_W(\ce)$. It is easily seen that 
$W(\ce,\theta)$ is isomorphic to a subgroup of $W_1$
(as all Cartan subalgebras
of $\z_\g(\ce)$ are $Z_G(\ce)$-conjugate), but these two groups
can be different in general.
We give below a sufficient condition for the equality to hold.

Let $\pi:\g\to \g\md G$ be the quotient mapping.
By \cite{ko63}, it is known that, for
$\xi\in\g\md G$, the fibre $\pi^{-1}(\xi)$ is an irreducible
normal complete intersection in $\g$ of codimension $l$. The complement of
the dense $G$-orbit in $\pi^{-1}(\xi)$ is of codimension at least 2.
Using a result of Richardson \cite{roger} and an extension of Springer's
theory to non-regular elements \cite{ls1}, 
we prove normality of some $G$-stable cones in $\g$ associated with 
$\theta$-groups.

\begin{s}{Theorem}  \label{norm3}
Suppose $\theta$ satisfies the relation $\dim\g_1\md G_0=k_{-1}$. 
Then, for any Cartan subspace $\ce\subset\g_1$, we have
\[
\pi^{-1}(\pi(\ce))=
\bigcap_{i:\ \esi_i\zeta^{d_i}\ne 1}\{x\in\g\mid F_i(x)=0\}  \ .
\]
This variety is irreducible, normal, and Cohen-Macaulay.
Furthermore, its ideal in $\bbk[\g]$ is generated by 
the above basic invariants $F_i$, i.e., $\pi^{-1}(\pi(\ce))$ is a complete 
intersection.
\end{s}\begin{proof}
Let $\te$ be a $\theta$-stable Cartan subalgebra containing $\ce$ and $W$
the corresponding Weyl group. Set $\sigma=\theta\vert_\te$.
By \cite[5.1]{ls1}, $W_1$ (which is not necessarily the same as either
$W^\sigma$ or $W(\ce,\theta)$)
is a reflection group in $\ce$ and the 
functions $F_i\vert_\ce$ with $\esi_i\zeta^{d_i}=1$
form a set of basic invariants for $W_1$
(our $k_{-1}$ is $a(d,\sigma)$ in \cite{ls1}). In other words,
$\bbk[\ce]^{W_1}$ is a graded polynomial algebra and the restriction mapping
$\bbk[\te]^W\to \bbk[\ce]^{W_1}$ is onto. 
This means that Theorem~B in \cite[\S\,5]{roger}
applies here, and we may conclude that $X:=\pi^{-1}(\pi(\ce))$
is normal and Cohen-Macaulay. Furthermore, Lemma~5.3 in loc.\,cit.
says that $X$ is irreducible; and the argument in p.~250 in loc.\,cit
shows that the ideal of $X$ is generated by the required basic invariants $F_i$.
\end{proof}%
We have proved before that the hypothesis of Theorem~\ref{norm3} is satisfied
for S-regular or \textrm{N}-regular gradings. However, in these cases some more
precise information is available.

\begin{s}{Theorem}  \label{norm}
Suppose $\theta\in\Aut\g$ is\/ \textup{\rmfamily N}-regular. Then
\begin{itemize}
\item[\sf (i)]\quad $G{\cdot}x\cap\g_1=G_0{\cdot}x$ for all $x\in\ce$. In
particular, $G{\cdot}x\cap\ce=G_0{\cdot}x\cap\ce$;
\item[\sf (ii)]\quad 
$W(\ce,\theta)=W_1$;
\item[\sf (iii)]\quad $\pi^{-1}(\pi(\ce))=
\ov{G{\cdot}\g_1}$ and all assertions of Theorem~\ref{norm3} hold
for this variety.
\end{itemize}
\end{s}\begin{proof}
There are two isomorphisms given by restriction

$\bbk[\g]^G \isom \bbk[\te]^{W}$ (Chevalley) and
$\bbk[\g_1]^{G_0}  \isom  \bbk[\ce]^{W(\ce,\theta)}$ (Vinberg, see 
\ref{vinberg}).
\\[.7ex]
Since $\theta$ is \textrm{N}-regular,         
$res_{\g,\g_1}: \bbk[\g]^G\to \bbk[\g_1]^{G_0}$ is onto by 
Theorem~\ref{onto+KW}(i). It follows that the restriction mapping
$res_{\te,\ce}: \bbk[\te]^W\to \bbk[\ce]^{W(\ce,\theta)}$ is onto, too.
In the geometric form, the ontoness of $res_{\g,\g_1}$ yields the closed
embedding $\g_1\md G_0\hookrightarrow \g\md G$. Because the points of such
(categorical) quotients parametrise the closed orbits and the closed 
$G_0$-orbits in $\g_1$ are those meeting $\ce$, the above embedding is
equivalent to the fact that $G{\cdot}x\cap\g_1=G_0{\cdot}x$ for all $x\in\ce$.
This gives (i). Similarly, the ontoness of $res_{\te,\ce}$ yields the
equality $W(\ce,\theta){\cdot}x=W\!{\cdot}x\cap\ce$ for all $x\in\ce$.
Since $W(\ce,\theta){\cdot}x\subset W_1{\cdot}x\subset W{\cdot}x\cap\ce$,
part (ii) follows.

Let us prove (iii). Set $X:=\pi^{-1}(\pi(\ce))$. It is a closed
$G$-stable cone in $\g$. By \cite[5.3]{roger}, $X$ is irreducible.
Since $\ov{G{\cdot}\g_1}$ is irreducible and 
$\ov{G{\cdot}\g_1}\subset X$, 
it suffices to verify that $\dim X=\dim\ov{G{\cdot}\g_1}$.
Because each fibre of $\pi$ is of dimension $\dim\g-l$ and 
$\dim\pi(\ce)=\dim\ce$, we obtain $\dim X=\dim\ce+\dim\g-l$. On the other hand,
regular elements are dense in $\g_1$, since $\theta$ is \textrm{N}-regular.
Therefore $\ov{G{\cdot}\g_1}$ contains a $\dim\ce$--parameter family 
of $G$-orbits of dimension $\dim\g-l$. This yields the required equality. 
In view of \ref{numerology}(v), Theorem~\ref{norm3} applies here.
\end{proof}%
{\it Remark.} For the S-regular locally free gradings, the coincidence of
$W(\ce,\theta)$ and $W_1$ was proved  in \cite[Prop.\,19]{vi}.
(In that case $W_1=W^\sigma$ for $\sigma=\theta\vert_{\z(\ce)}$.)
Therefore, in view of Theorem~\ref{main}(iii), the equality
\ref{norm}(ii) is an extension of that result of Vinberg.

\begin{s}{Proposition}  \label{norm2}
Suppose $\theta$ is\/ \textup{\rmfamily S}-regular. Let $\ce$ be any Cartan subspace in 
$\g_1$. Then 
\[
\pi^{-1}(\pi(\ce))=\ov{G{\cdot}\g_1}= \ov{G{\cdot}\ce} \ 
\]
and all assertions of Theorem~\ref{norm3} hold.
\end{s}\begin{proof*}
The first equality stems from the presence of regular elements in $\g_1$
(cf. the proof of \ref{norm}(iii)). Since regular semisimple elements are
dense in $\g_1$, $\ov{G_0{\cdot}\ce}=\g_1$. This gives the second equality.
In view of \ref{main}(i), Theorem~\ref{norm3} applies here.
\end{proof*}%
%
\begin{rem}{Examples}
1. Consider again the automorphism $\theta_2$ from
Example~\ref{primer}. As we already know, 
$\theta_2$ is not \textrm{N}-regular and the ontoness of 
$res_{\g,\g_1}$ fails here. The latter shows that 
 \ref{norm}(i)  does not hold.
Since $\ce$ contains regular elements, $Z_W(\ce)=\{1\}$. It is easily seen
that $W_1$ is isomorphic to $W^{\sigma}$, the centraliser of 
$\sigma=\theta\vert_\te$ in $W$. Springer's theory \cite[\S\,4]{tony}
says that $W^{\sigma}\subset GL(\ce)$ is a finite reflection group whose
degrees are those degrees of $W$ that are divisible by $m$, i.e., 8,\,12.
Hence $\#W_1=96$ and $\#W(\ce,\theta)=48$. Thus, \ref{norm}(ii) fails, too.
However, $\pi^{-1}(\pi(\ce))$ is normal, Cohen-Macaulay, etc.,
and the reason is that $\theta_2$ is
S-regular.

2. It really may happen that $\theta$ is neither \textrm{N}-regular nor S-regular, but
the equality $\dim\g_1\md G_0=k_{-1}$ holds. Let $\g$ be a simple Lie algebra
of type ${\Bbb E}_7$. Consider the inner automorphism $\theta$ that is determined
by the following Kac's diagram: \\
\centerline{
\raisebox{12\unitlength}{$\theta$:}\quad 
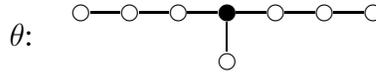
\begin{picture}(90,40)(0,0)
\setlength{\unitlength}{0.017in}
\put(50,8){\line(0,1){9}}
\put(50,20){\circle*{5}} 
\multiput(8,20)(15,0){6}{\line(1,0){9}}
\multiput(5,20)(15,0){7}{\circle{5}}
\put(50,5){\circle{5}}
\end{picture} 
}

\noindent Then the order of $\theta$ is 4, $G_0=A_3\times A_3\times A_1$, 
and $\dim\g_1=32$. Here $k_{-1}=\dim\g_1\md G_0=2$, so that 
Theorem~\ref{norm3} applies. On the other hand, $(\dim\g-\rk\g)/4\notin
{\Bbb N}$. Hence $\theta$ is not S-regular. It can be shown that the maximal
nilpotent orbit meeting $\g_1$ is of dimension 120 (its 
Dynkin-Bala-Carter label is $E_6$), i.e., $\theta$ is not \textrm{N}-regular.
\end{rem}
\vskip-1ex


\section{Exponents and coexponents of little Weyl groups}
\setcounter{equation}{0}

\noindent
We keep the previous notation. 
Here we briefly discuss some other consequences of \cite{ls1},\,\cite{ls2}
for $\theta$-groups.
\\[.7ex]
Recall the definition of (co)exponents.
If $\tilde W$ is a reflection group in $\ce$, then 
$(\bbk[\ce]\otimes N)^{\tilde W}$ is a graded free $\bbk[\ce]^{\tilde W}$-module
for any $\tilde W$-module $N$. The {\it exponents} (resp. {\it coexponents})
of $\tilde W$ are
the degrees of a set of free homogeneous generators for
this module, if $N=\ce^*$ (resp. $N=\ce$). 
As is well known, if $\{\tilde d_i\}$ are the degrees of basic invariants
in $\bbk[\ce]^{\tilde W}$, then $\{\tilde d_i-1\}$ are the exponents.

The theory of Lehrer and Springer gives a description of coexponents for
the subquotient $W_1=N_W(\ce)/Z_W(\ce)$ under the constraint
$\dim\ce=k_{-1}$, while the theory of $\theta$-groups deals with the
group $W(\ce,\theta)=N_{G_0}(\ce)/Z_{G_0}(\ce)$.
By \ref{numerology}(v) and \ref{norm}(ii), we know that $\dim\ce=k_{-1}$ and
$W_1=W(\ce,\theta)$ whenever $\theta$ is \textrm{N}-regular.
Thus, \textrm{N}-regularity allows us to exploit the theory of Lehrer and
Springer in the context of 
$\theta$-groups. However, to use that theory in full strength,
we need the constraint that $\theta$ is S-regular, too.
\begin{s}{Proposition}
Suppose $\theta\in\Aut\g$ is\/ \textup{\rmfamily N}-regular and let $\{e,h,f\}$
be a $\theta$-adapted regular $\tri$-triple.
\begin{itemize}
\item[\sf (i)]\quad
The exponents of $W(\ce,\theta)$ correspond to the eigenvalues of $\theta$
on $\z(e)_{-1}$. More precisely, $m_j$ is an exponent of $W(\ce,\theta)$
if and only if $\esi_j\zeta^{m_j}=\zeta^{-1}$;
\item[\sf (ii)]\quad If $\theta$ is also S-regular, then
the coexponents of $W(\ce,\theta)$ correspond to the eigenvalues of $\theta$
on $\z(e)_{1}$. More precisely, $m_j$ is a coexponent of $W(\ce,\theta)$
if and only if $\esi_j\zeta^{m_j}=\zeta$.
\end{itemize}
\end{s}\begin{proof}
(i) This part is essentially contained in Theorem~\ref{onto+KW}(i). \\
(ii) A description of coexponents for subfactors of the form 
$N_W(\ce)/Z_W(\ce)$, if $\ce$ contains a regular vector,
is due to Lehrer and Springer. However, an explicit formulation was given
only in the ``untwisted'' case (see Theorem C in \cite{ls2}),
which in our setting correspond to the case where $\theta$ is inner.
A general statement can be extracted from the discussion
in \cite[\S\,4]{ls2}, especially from Propositions 4.6, 4.8, and
Theorem D. 

The connection with $\theta$-eigenvalues on $\z(e)$ follows from 
Theorem~\ref{numerology}. 
\end{proof}%
{\it Remark.}
If $\theta$ is \textrm{N}-regular, then formulas of Section~\ref{section-R}
shows that
$k_{-1}=\dim\ce=\dim\z(e)_{-1}$ and
$k_1=\dim\z(e)_1$.
S-regularity guarantee us that $k_{-1}=k_1$ and, even stronger, that
$k_{-i}=k_i$ for all $i$ (see \re{main}).
However it can happen that $k_1> k_{-1}$. For instance, if $\theta$ is an
\textrm{N}-regular inner automorphism of ${\Bbb E}_6$, of order 5, then $k_1=2$
and $k_{-1}=1$.
In this case, it is not clear how
to characterise the coexponents for $W(\ce,\theta)$.
Another difficulty may occur if $\theta$ is S-regular, but not \textrm{N}-regular.
Here $k_{-1}=k_1$, but the groups $W_1$ and $W(\ce,\theta)$ can be
different.

\begin{rem}{Example}
The exceptional Lie algebra $\g={\Bbb E}_{6}$ has an outer automorphism 
$\theta$ of order 4 such that $G_0=A_3\times A_1$. 
This automorphism is determined by the following Kac's diagram
(the underlying graph is the Dynkin diagram of type $E_6^{(2)}$)\ : \quad
\begin{picture}(100,25)(0,5)
\multiput(10,8)(20,0){5}{\circle{6}}
\put(70,8){\circle*{6}}
\multiput(13,8)(20,0){2}{\line(1,0){14}}
\put(73,8){\line(1,0){14}}
\multiput(52.5,7)(0,2){2}{\line(1,0){15}}
\put(55,5){$<$} 
\end{picture}.
\vskip1ex\noindent
Here $\g_1$ is the tensor product of a 10-dimensional representation of
$A_3$ (with highest weight $2\varphi_1$) and the simplest representation 
of $A_1$.
We have $\dim\g_0=18$ and $\dim\g_1=20$.
It is easily seen that the $G_0$-representation on $\g_1$ is locally free,
hence $\dim\g_1\md G_0=\dim\ce=2$. Using Proposition~\ref{1/m}(i), we conclude
that $\theta$ is S-regular and then, by Theorem~\ref{main}(iii), that
$\theta$ is also \textrm{N}-regular. In order to use the preceding Proposition,
one has to know the factors $\{\esi_i\}$. In this case the pairs
$(m_i,\esi_i)$ ($1\le i\le 6$) are  $(1,1),(4,-1),(5,1),(7,1),(8,-1),(11,1)$.
Then an easy calculation shows that the exponents of $W(\ce,\theta)$ are
7,\,11 and the coexponents are 1,\,5. Then
looking through the list of the irreducible finite reflection groups,
one finds that here $W(\ce,\theta)$ is Group 8 in the Shephard-Todd numbering.
(A list including both the exponents and the coexponents is found in 
\cite[Table 2]{os}).  
\end{rem}%
%

%
%

\begin{thebibliography}{LS1}

\bibitem[An]{leva2} {\sc L.V.~Antonyan}.  On homogeneous elements of
periodically graded semisimple Lie algebras, in: 
``Questions in group theory and homological algebra'' (A.L.~Onishchik, Ed.),
pp. 55--64. Yaroslavl, 1987 (Russian).



\bibitem[Ko]{ko63} {\sc B.~Kostant}. Lie group representations in
polynomial rings, {\it  Amer. J. Math.} {\bf 85}(1963), 327--404.

\bibitem[KR]{KR} {\sc B.~Kostant, S.~Rallis}. Orbits and 
representations associated with symmetric spaces, 
{\it Amer. J. Math.} {\bf 93}(1971), 753--809.

\bibitem[LS1]{ls1} {\sc G.~Lehrer, T.A.~Springer}.
Intersection multiplicities and reflection subquotients of unitary 
reflection groups,~I. {\it In}: ``Geometric group theory down under'' 
(J.~Cossey, Ed.) W. de Gruyter (1999), 181--193.

\bibitem[LS2]{ls2} 
{\sc G.~Lehrer, T.A.~Springer}.
Reflection Subquotients of Unitary Reflection Groups,
{\it Can. J. Math.}  {\bf 51}(1999), 1175--1193.

\bibitem[OS]{os}  {\sc P.~Orlik, L.~Solomon}.
Unitary Reflection Groups and Cohomology, {\it Invent. Math.}
{\bf 59}(1980), 77--94.
 

\bibitem[Pa]{R=S} {\sc D.~Panyushev}. 
Regular elements in spaces of linear representations II,
{\it Izv. AN SSSR, Ser. matem.}
{\bf 49}(1985), {\rus N0}\,5, 979--985 (Russian). 
English translation: {\it Math. USSR-Izv.} {\bf 27}(1986), 279--284.


\bibitem[Po]{cov} {\sc V.L.~Popov}. Representations with a free module of
covariants, {\it Funkt. Analiz i ego Priloz.} {\bf 10}:3 (1976), 91--92
(Russian). English translation: {\it Funct. Anal. Appl.} {\bf 10}(1976),
242--245. 

\bibitem[Ri]{roger} {\sc R.W.~Richardson}.
Normality of $G$-stable subvarieties of a semisimple Lie algebra,
In: {\it ``Algebraic Groups''},  
Lecture Notes Math., Springer-Verlag {\bf 1271}(1987), 243--264.

\bibitem[Sp]{tony} {\sc T.A.~Springer}.
Regular elements of finite reflection groups, {\it  Invent. Math.}
{\bf 25}(1974), 159--198.


\bibitem[St]{endo} {\sc R.~Steinberg}.  
Endomorphisms of algebraic groups, {\it Memoirs A.M.S.}
{\bf 80}(1968), 1--108 (= {\it Collected Papers}, A.M.S., 1997, pp.~229--286).

\bibitem[Vi1]{vi}
{\sc E.B.~Vinberg}. The Weyl group of a graded Lie algebra,
{\it Izv. AN SSSR. Ser. Matem.} {\bf 40}(1976), {\rus N0}\,3, 488--526
(Russian). English translation: {\it Math. USSR-Izv.} {\bf 10}(1976), 463--495.

\bibitem[Vi2]{vi79}
{\sc E.B.~Vinberg}. Classification of homogeneous nilpotent elements 
of a semisimple graded Lie algebra, In: {\it ``Trudy seminara po vect. 
i tenz. analizu"}, vol.\,19, pp.~155--177. Moscow: MGU 1979 (Russian). 
English translation:
{\it Selecta Math. Sovietica} {\bf 6}(1987), 15--35.

\bibitem[VO]{vion}
{\sc E.B.\,Vinberg} and {\sc A.L.\,Onishchik}.
``Seminar on Lie groups and 
algebraic groups",  Moskva: ``Nauka'' 1988 (Russian).
English translation: {\sc A.L.\,Onishchik} and {\sc E.B.\,Vinberg}:
``Lie groups and algebraic groups''. Berlin Heidelberg 
New York: Springer 1990.

\bibitem[VP]{vipo}
{\sc E.B.~Vinberg, V.L.~Popov}. ``Invariant theory", In: {\it 
Sovremennye problemy matematiki. Fundamentalnye napravleniya}, t.\,55,
pp.\,137--309. Moscow: VINITI 1989 (Russian).
English translation in: Algebraic Geometry IV
(Encyclopaedia Math. Sci., vol.~55, pp.123--284) 
Berlin Heidelberg New York: Springer 1994.

\end{thebibliography}
\end{document}